\documentclass[preprint,12pt]{elsarticle}

\usepackage{amssymb}
\usepackage{amsmath}

\usepackage[american]{babel}

\usepackage{natbib} 

\usepackage{mathtools} 
 \usepackage{siunitx} 
\mathtoolsset{showonlyrefs}
\usepackage{booktabs} 
\usepackage{tikz} 
\usepackage{amsthm}
\usepackage{txfonts}
\newcommand{\nodeset}{\mathcal{N}}
\newcommand{\nn}[0]{\hspace*{.7em}}

\newtheorem{prop}{Proposition}
\newtheorem{coro}{Corollary}

\newtheorem{thm}{Theorem}
\theoremstyle{remark}
\newtheorem{example}{Example}

\journal{Statistics \& Probability Letters}

\begin{document}

\begin{frontmatter}

\title{On Exponential Random Graph Models with Dyadic Independence} 

\author{Kayvan Sadeghi} 

\affiliation{organization={Department of Statistical Science, University College London},
            addressline={Gower St}, 
            city={London},
            postcode={WC1E 6BT}, 
            country={United Kingdom}}

\begin{abstract}
We show that the only exponential random graph model with $n$ nodal parameters, dyads being independent,  and the natural assumption of permutation-equivariant nodal parametrization  is the $\beta$ model. In addition, we show that an exponential random graph model with similar assumptions but with fewer than $n$ block parameters is the additive stochastic block model. We also provide similar results for directed networks.
\end{abstract}




\end{frontmatter}



%
%
%

\section{Introduction}
Exponential random graph models (ERGMs) \citep{frank:91,wasserman:pattison:96} are some of the most used statistical models for network data. These are statistical models in exponential family form with sufficient statistics that are some feature of the observed network(s).

The  $\beta$ model  \citep{bli10,chatterjee:diaconis:sly:11} is one of the most used types of ERGMs, where there are $n$ parameters corresponding to the nodes of the network and the sufficient statistics for the model is the degree sequence of the network.  Stochastic block models \citep{SBM:83} are another important type of ERGMs, where parameters are assigned to blocks of nodes. The parameters associated to the nodes in these models are called ``nodal'' parameters, as they express the propensity of the nodes---as in the $\beta$ model case---, and  the parameters associated to the blocks may be called ``block'' parameters---as in the stochastic block model case. In the latter case, one can think of the block parameters as nodal parameters with some additional constrains that set the nodal parameters to be equal in each block.

When the ERGM deals with directed random networks as opposed to undirected ones, the mentioned models are easily generalizable. For example, the directed version of the $\beta$ model has been known and used for a long time under the name of the $p_1$ model \citep{hol81,fie81}.

An important property of these mentioned models, as opposed to most ERGMs, is that they are ``dyadic independent'', meaning that, in these models, each edge existing is mutually independent on the existence of other edges in the network.

Some of the popularity of these particular ERGMs is due to their desirable properties: The existence of the MLE is well-understood for these models \citep{rin13}, and there are simple fixed-point algorithms to find the MLE \citep{chatterjee:diaconis:sly:11}; the asymptotics of these models are known \citep{che21,cel12}, and they are generalizable to the higher order random hyper-graphs \citep{sta15,gho14}. Although well-behaved and well-used, there is little justification in the literature as why these models, and the strong modeling assumptions behind them, might be reasonable or justified to assume.

In this paper, we attempt to address this issue. We use simple and elementary derivations to prove that the only dyadic independent ERGMs with $n$  (i.e., the number of nodes) nodal parameters and with certain natural assumption of ``permutation-equivariance''  must be the $\beta$ model. In addition, the only dyadic-independent stochastic block models with the assumption of permutation-equivariance with fewer than $n$ parameters (corresponding to the blocks) are the additive stochastic block model. Additive stochastic block models resemble the $\beta$ model with nodes being replaced by blocks of nodes. 

We rather trivially extend the results to the case of random directed networks, where under the same assumptions we obtain the $p_1$ model and a directed version of the additive stochastic block model.

The assumption of permutation-equivariance for nodal parametrization is quite natural in this setting. It implies that the type of relations between the dyadic parameters and the corresponding nodal parameters must be the same across the network, so although the values of the paprameters can be different they are all mapped to the nodal parameters in a similar way. Although, equivariance has been used in statistics for estimators \citep{leh05}, and is becoming an increasingly popular assumption for neural networks \citep{vic21}, as far as we know, they have not been applied to parameters as opposed to the data. We believe, such an assumption on parameters could be interesting in its own right in settings beyond random network analysis.

Although the results are quite easy to prove, and seem rather basic, as discussed, they are important for justifying the use of $\beta$ models and $p_1$ models as well as additive stochastic block models. To the best of our knowledge, these results have not been known, or at least explicitly stated, in the literature.
   
The structure of the paper is as follows: In the next section, we provide an introduction to the definitions and results needed on dyadic-independent exponential families, $\beta$- and $p_1$- models as well as the stochastic block model. In Section \ref{sec:res}, we provide the definition of permutation-equivariance for the nodal parametrization, and provide and prove our results. We conclude the paper in Section \ref{sec:sum}.
\section{Dyadic-independent exponential random graph models}
\subsection{Exponential random graph models}\label{sec:5}

Given a finite or countably infinite node set $\nodeset$ ---
representing \emph{individuals} or \emph{actors} in a given population of
interest ---  we define a
\emph{random network} over $\nodeset$ to be a collection $X = (X_d, d\in
D(\nodeset) )$
of binary random variables taking values $0$ and $1$ indexed by a set
$D(\nodeset)$ of \emph{dyads}. Throughout the paper we let $|\nodeset|=n$, and take
$D(\nodeset)$ to be the set of all
pairs $ij$ of nodes in $\nodeset$, where $i$ and $j$ are distinct, and
nodes $i$ and $j$ are said to have
an \emph{edge}  if the random variable $X_{ij}$ takes the value $1$, and no edge
otherwise. For \emph{directed networks}, i.e., when the edges are assumed to be directed, the $ij$ pair is ordered as $(i,j)$, indicating the direction of the edge from $i$ to $j$. When dealing with \emph{undirected networks}, i.e., when the edges are assumed to be undirected, we let the $ij$ pair to be unordered.

Thus, in general, a network is a random
variable taking value in $\{0,1 \}^{n(n-1)}$ and can therefore
be seen as a random {\it labelled} graph with node set $\nodeset$,
whereby the edges form the random edges of the graphs.
We will write $\mathcal{G}_{\nodeset}$ for the set of all simple, labelled graphs on
$\nodeset$. As discussed, for a network $x\in\mathcal{G}_{\nodeset}$, we allow in general to have $x_{ij}\neq x_{ji}$, but for undirected networks, we impose $x_{ij}= x_{ji}$.




\emph{Exponential random graph models} \citep{frank:91,wasserman:pattison:96},
or ERGMs in short, are among the most important and popular statistical models
for network data.

ERGMs are exponential families of distributions \citep{barndorff:78} on
$\mathcal{G}_{\nodeset}$  whereby the probability of observing a network $x \in
\mathcal{G}_{\nodeset}$ can
be expressed as
\begin{equation}
    P_{\theta}(x) = \exp\left\{\sum_{l=1}^ms_l(x)\theta_l-\psi(\theta)\right\},
    \quad \theta \in \Theta \subseteq \mathbb{R}^m,\label{eq:ergm}
\end{equation}
where $s(x) = (s_1(x) \ldots, s_m(x)) \in \mathcal{R}^m$ are \emph{canonical sufficient
statistics} which capture some important feature of $x$, $\theta \in
\mathbb{R}^m$ is a point in the canonical parameter space $\Theta$ and
$\psi \colon \Theta \rightarrow [0,\infty)$ is the \emph{log-partition function}, also known as the \emph{normalizing constant}.

The choice of the canonical sufficient statistics and any restriction imposed on
$\Theta$ will determine the properties of the corresponding ERGM.

\subsection{Dyadic Independent ERGMs}
We call an ERGM \emph{dyadic independent} if the dyads $(X_{ij}, i,j\in\mathcal{N})$ are mutually independent. Denote by $p_{ij}$ the probability that the dyad $X_{ij}$ takes the value $1$, i.e., $p_{ij}=P(X_{ij}=1)$. First, for an undirected network $x \in
\mathcal{G}_{\nodeset}$,  from dyadic independence, we have:

\begin{align}
P(x)&=\prod_{i,j:i\neq j}p_{ij}^{x_{ij}}(1-p_{ij})^{x_{ij}}=\prod_{i,j:i\neq j}(\frac{p_{ij}}{1-p_{ij}})^{x_{ij}}(1-p_{ij})\\
&=\exp\left\{\sum_{i,j:i\neq j}x_{ij}\log\frac{p_{ij}}{1-p_{ij}}+\log(1-p_{ij})\right\}=\exp\left\{\sum_{i,j:i\neq j}x_{ij}\theta_{ij}-\psi(\theta)\right\},\label{eq:satur}
\end{align}
where  $\theta_{ij}=\log p_{ij}/(1-p_{ij})$ are $n(n-1)$ canonical parameters, known as the \emph{logits}, and $\psi(\theta)$ is a normalizing constant. We see that for this saturated ERGM with dyadic independence, the minimal sufficient statistics are all the dyads.

In the case of undirected networks, in \eqref{eq:satur}, we have that $p_{ij}=p_{ji}$. This is the same as the products and sums in \eqref{eq:satur} being over $\{i,j:i<j\}$. Hence, in this case, we end up with the additional constraint that $\theta_{ij}=\theta_{ji}$, and the parameter space ends up to be  ${n \choose 2}$ dimensional.
%
\subsection{Some dyadic-independent ERGMs}
The simplest ERGM, for undirected networks, is the \emph{Erd\"{o}s--R\'{e}nyi} model \citep{erdos:renyi:60}, where there is only one parameter $\theta$, and the canonical sufficient statistic is the number of edges in $x$. This is equivalent to edges occurring independently with probability $p$ and in Equation \eqref{eq:satur}, for all $i,j\in \nodeset$, $\theta_{ij}=\theta$.

\emph{Beta models} \citep{chatterjee:diaconis:sly:11,rin13} (also written $\beta$ models) are also undirected ERGMs, and they can be considered a generalization of the Erd\"{o}s-R\'{e}nyi model: in these models it is also assumed that edges occur independently; the canonical parameters $\theta_{ij}$ in Equation \eqref{eq:satur} are given as $\theta_{ij}=\beta_i+\beta_j$,
where  $(\beta_i, i\in\nodeset)$ can be interpreted as  parameters that
determine the propensity of node $i$ to have edges. The probability  of a
network $x \in \mathcal{G}_{\nodeset}$ is thus
\begin{align}
P(x)&=\exp\left\{\sum_{i,j:i<j}x_{ij}(\beta_i+\beta_j)-\psi(\beta)\right\}=
\exp\left\{\sum_{i\in\nodeset}d_i(x)\beta_i-\psi(\beta)\right\},\label{eq:betarep}
\end{align}
where $(d_i(x), i\in \nodeset)$ is  the \emph{degree sequence} of $x$, i.e.\ $d_i(x)$ is the number of edges in $x$ involving node $i$.

The \emph{stochastic block model} (SBM) \citep{SBM:83} (for undirected networks) is defined as follows:   In  a SBM, the nodes are grouped in $r$ groups, called \emph{blocks} or \emph{communities}. The edges are then occurring independently with the same probability of occurring within each block and the same probability of occurring between each pair of blocks.

More formally, we say that group labels are given by map $b:\nodeset\rightarrow [r]$, where $[r]=\{1,\cdots,r\}$. The fibre over  $l\in[r]$ is a block of nodes, which we denote by $z_l$. We let probabilities of dyads being $1$ be $p_{ij}= q_{b(i)b(j)}=q_{kl}$ for block labels $k,l$, to which $i,j$ belong respectively. This implies that  $\theta_{ij}= \eta_{b(i)b(j)}:=\eta_{kl}$, where $\eta_{kl}=\log q_{kl}/(1-q_{kl})$. From  \eqref{eq:satur}, we obtain
\begin{align}
P(x)&= \exp\left\{\sum_{k,l:k\leq l}\sum_{i,j: i\in z_k,j\in z_l}x_{ij}\eta_{kl}-\psi(\eta)\right\}
= \exp\left\{\sum_{k,l:k\leq l} e_{kl}(x)\eta_{kl}-\psi(\eta)\right\},\label{eq:sbm}
\end{align}
where $e_{kl}(x)$, for $k\neq l$, is the number of edges between the blocks $z_k$ and $z_l$, and  $e_{kk}(x):=e_k(x)$ is the number of edges within $z_k$.  This is the form of stochastic block model we shall use.

Now, in the same fashion as in the $\beta$ model case, if there is a further parametrization $\eta_{kl}=\delta_k+\delta_l$, for $k\neq l$ then
\begin{align}
P(x)&=\exp\left\{\sum_{k,l:k< l}e_{kl}(\delta_k+\delta_l)+e_k(x)\eta_{kk}-\psi(\delta,\eta)\right\}\\ & =
\exp\left\{\sum_k d_k(x)\delta_k+e_k(x)\eta_{kk}-\psi(\delta,\eta)\right\},\label{eq:addsbm}
\end{align}
where $d_k(x)$ is the degree of block $z_k$. This gives the block version of $\beta$ model called the \emph{additive stochastic block model}.


The directed version of the $\beta$ model is the \emph{configuration model} which is a submodel of the $p_1$ model of \citet{hol81}; see also \citet{fie81} for a modification. In this model $\theta_{ij}=\alpha_i+\beta_j$, which leads to
\begin{align}
P(x)&=\exp\left\{\sum_{i,j: i\neq j}x_{ij}(\alpha_i+\beta_j)-\psi(a,b)\right\} =
\exp\left\{\sum_{i\in\nodeset}a_i(x)\alpha_i+b_i(x)\beta_i-\psi(\alpha,\beta)\right\},\label{eq:betarepd}
\end{align}
where $a(x)$ is  the \emph{out-degree sequence} of $x$, i.e.\ $a_i(x)$ is the number of edges in $x$ pointing out from node $i$; and $b(x)$ is  the \emph{in-degree sequence} of $x$, i.e.\ $b_i(x)$ is the number of edges in $x$ pointing into node $i$. The difference of this model with the $p_1$ model is that it does not contain a \emph{reciprocal parameter} with corresponding sufficient statistics $a(x)b(x)$, which make dyads $x_{ij}$ and $x_{ji}$ dependent. In addition, in the $p_1$ model $\alpha$ and $\beta$ are often written as the sum of two parameters (e.g. $\mu+\alpha'$ in the former case), where $\mu$  presents the global density of the network.

Following the same procedure, we define a \emph{directed additive stochastic block model} for a directed network to be
\begin{equation}
P(x)=\exp\{\sum_{k\in [r]}[a_k(x)\delta_k+b_k(x)\lambda_k+e_k(x)\eta_{kk}-\psi(\delta,\lambda,\eta)]\},\label{eq:addsbmd}
\end{equation}
where $a_k(x)$ and $b_k(x)$ are the out-degrees and in-degrees of the block $z_k$.

\section{Permutation-equivariant dyadic-independent exponential random graph models}\label{sec:res}
\subsection{Permutation-equivariant nodal parametrization}
Assume that a group $G$ acts (from the left) on the set $X$. Then a function $f : X \rightarrow Y$ is said to be  \emph{invariant} under $G$ if $f(\pi\cdot x) = f(x)$ for all $\pi\in G$ and all $x\in X$.

Assume also that the group $G$ acts (from the left) on the set $Y$. Then $f$ is said to be \emph{equivariant} under $G$ if $f(\pi\cdot x) = \pi\cdot f(x)$ for all $\pi\in G$ and all $x\in X$.

In this paper, we concern an extension of equivariance to arrays, and particularly, with permutation groups $\pi:\mathcal{N}\rightarrow\mathcal{N}$.
Consider a parametrization $\theta_{B}=f_{(B)}(\delta_A)$, where $\delta_A=(\delta_j)_{j\in A}$, for $A\subset \nodeset$, $i\in\nodeset$, and similarly for $\theta$, and $\theta\in \Theta$ and $\delta\in \Delta$. For a permutation $\pi:\nodeset\rightarrow\nodeset$, we define $\pi(\delta_A)=\delta_{\pi(A)}$; and similarly $\pi(\theta_B)=\theta_{\pi(B)}$; hence, we have that $\pi(f_{(B)}(\delta_A))=\theta_{\pi(B)}$. We say that this parametrization  is \emph{permutation equivariant} if for every permutation $\pi$, all $f_{(B)}$ are permutation equivariant, i.e., $\pi(f_{(B)}(\delta_A))=f_{(B)}(\pi(\delta_A))$.

In the particular case of dyadic-independent ERGMs, from Equation \eqref{eq:satur}, it is evident that a further parametrization of the dyadic parameters to $n$ nodal parameters $(\theta_i)_{i\in \nodeset}$ is of form $\theta_{ij}=f_{(ij)}(\theta_i,\theta_j)$, $i,j\in\nodeset$.   Similar to the previous case, for a permutation $\pi$, we define $\pi(\theta_i,\theta_j)=(\theta_{\pi(i)},\theta_{\pi(j)})$; and similarly $\pi(\theta_{ij})=\theta_{\pi(i)\pi(j)}$. We now say that this parametrization is \emph{(doubly) permutation equivariant} if for every permutation $\pi$, all $f_{(ij)}$ are (doubly) permutation equivariant, i.e., $\pi(f_{(ij)}(\theta_i,\theta_j))=f_{(ij)}(\pi(\theta_i,\theta_j))$.

The assumption of permutation-equivariance in this case leads to the equality of all the functions $f_{(ij)}$:
\begin{prop}\label{samef}
If the  parametrization $\theta_{ij}=f_{(ij)}(\theta_i,\theta_j)$, $i,j\in\nodeset$, from dyadic to nodal parameters in a dyadic-independent ERGM is permutation equivariant, then for all $i,j$, $f_{(ij)}\equiv f$, for a fixed function $f:\mathbb{R}^2\rightarrow\mathbb{R}$.
\end{prop}
\begin{proof}
We have that
$$f_{(ij)}(\pi(\theta_i,\theta_j))=\pi(f_{(ij)}(\theta_i,\theta_j))=\pi(\theta_{ij}).$$
For arbitrary $i,j,k,l$, if $\pi$ maps $(i,j)$ to $(k,l)$ then we have
$$f_{(ij)}(\theta_k,\theta_l))=\theta_{kl}=f_{(kl)}(\theta_k,\theta_l).$$
We denote all these functions by $f$.
\end{proof}

\subsection{The case of $n$ nodal parameters: $\beta$ model}
As above, we assume that there are $n$ nodal parameters $(\theta_i)_{i\in \nodeset}$, and that the ERGM is dyadic independent with permutation-equivariant nodal parametrization. We further assume that the ERGM is defined on undirected networks. 
\begin{thm}\label{thm:beta}
Assume an ERGM for undirected networks with $n$ nodal parameters is dyadic independent and obtained from a permutation-equivariant parametrization $\theta_{ij}=f_{(ij)}(\theta_i,\theta_j)$, $i,j\in\nodeset$, where $f_{(ij)}$ are differentiable. It then holds that $\theta_{ij}=g(\theta_i)+g(\theta_j)$, for all $i,j\in\nodeset$ and some function $g$.
\end{thm}
\begin{proof}
By Proposition \ref{samef}, $\theta_{ij}=f(\theta_i,\theta_j)$.
By minimal sufficiency theorem and Equation \eqref{eq:satur}, we have that, for every $x,y$,
$$\exp\left\{\sum_{i,j:i\neq j}f(\theta_i,\theta_j)\cdot(x_{ij}-y_{ij})\right\}=c,$$
where $c$ is a constant in $\theta$. Hence, by differentiation over an arbitrary $\theta_k$,
$$\sum_{j: j\neq k} df(\theta_k,\theta_j)/d\theta_k \cdot(x_{kj}-y_{kj})=0,$$
for every values of $\theta_k$ and $\theta_j$. This implies that for every $j$, $df(\theta_k,\theta_j)/d\theta_k$ is constant in $\theta_j$. Therefore, $f(\theta_i,\theta_j)=h(\theta_i)+g(\theta_j)$ for some functions $g,h$.

Now since $\theta_{ij}=\theta_{ji}$, we have that,  $h(\theta_i)+g(\theta_j)=h(\theta_j)+g(\theta_i)$, for all $i,j\in\nodeset$ and for every values of $\theta_i$ and $\theta_j$. This implies that $h\equiv g$.
\end{proof}
\begin{coro}\label{coro:beta}
A dyadic independent ERGM for undirected networks with $n$ nodal parameters obtained from a permutation-equivariant parametrization is a $\beta$ model.
\end{coro}
\begin{proof}
Since $f$ is smooth, $g$ has to be smooth. Hence, we simply define $\beta_i=g(\theta_i)$. Now  \eqref{eq:betarep}, gives the $\beta$ model. 
\end{proof}
The following example illustrates that permutation-equivariance is necessary:
\begin{example}
Let $\nodeset=\{1,\cdots,n\}$. Consider a further parametrization of $\theta_{ij}$ in a dyadic-independent ERGM as follows:
$$\theta_{1j}=\theta_j,\nn j>1;\nn \theta_{ij}=\theta_1\nn i>1, \forall j>i.$$
This parametrization is nodal as each  $\theta_{ij}$ depends on the parameters associated with the endpoints (with some additional equality constraints). This model is clearly different from the $\beta$ model. Notice also that this parametrization is not permutation equivariant: for a permutation $\pi$ that maps $1$ to $2$, and $2$ to $3$, $\theta_{12}=\theta_2$, but $\pi(\theta_{12})=\theta_1\neq \pi(\theta_2)=\theta_3$.
\end{example}
\subsection{The case of fewer than $n$ block parameters: additive stochastic block model}
In the same manner as nodal parameters, we can consider fewer than $n$ block parameters to be a further parametrization of the parameters $(\eta_{kl})_{k,l\in [r]}$ in the stochastic block model of \eqref{eq:sbm}.

In particular, in the case of dyadic-independent stochastic block models, further parametrization to $r<n$ block parameters is of form $\eta_{kl}=h_{(kl)}(\eta_k,\eta_l)$, $k,l\in[r]$.   Similar to the previous case, for a permutation $\pi$, we define $\pi(\eta_k,\eta_l)=(\eta_{\pi(k)},\eta_{\pi(l)})$; and similarly $\pi(\eta_{kl})=\eta_{\pi(k)\pi(l)}$. This parametrization is now (doubly) permutation equivariant if for every permutation $\pi$, all $h_{(kl)}$ are (doubly) permutation equivariant, i.e., $\pi(h_{(kl)}(\eta_k,\eta_l))=h_{(kl)}(\pi(\eta_k,\eta_l))$.

Now, we can prove in a similar fashion to Theorem \ref{thm:beta}, the following theorem:
\begin{thm}\label{thm:sbm}
If an undirected dyadic-independent stochastic block model with $r<n$ block parameters is obtained from a permutation-equivariant parametrization $\eta_{kl}=h_{(kl)}(\eta_k,\eta_l)$, $k,l\in[r]$, where  $h_{(kl)}$ are differentiable, then $\eta_{kl}=g(\eta_k)+g(\eta_l)$, for all $k,l\in [r]$ and some function $g$.
\end{thm}
And similar to Corollary \ref{coro:beta}, we have
\begin{coro}
A dyadic independent stochastic block model with $r<n$ block parameters obtained from a permutation-equivariant parametrization is an additive stochastic block model.
\end{coro}
\subsection{The cases of directed networks}
The previous case of undirected networks naturally generalizes to the case of directed networks. We note that the definition of doubly permutation equivariance is independent of the type of the network, and Proposition \ref{samef} holds regardless of whether the network is undirected or directed. Here, we provide the results presented previously for directed networks. We refrain from providing the proofs as they are a trivial extension of the proofs for the undirected case. 

For the case of $n$ nodal parameters $(\theta_i)_{i\in \nodeset}$, and the ERGM that is dyadic independent with permutation-equivariant nodal parametrization, but this time for directed netwroks, we have the following:
\begin{thm}\label{thm:p1}
Assume an ERGM for directed networks with $n$ nodal parameters is dyadic independent and obtained from a permutation-equivariant parametrization $\theta_{ij}=f_{(ij)}(\theta_i,\theta_j)$, $i,j\in\nodeset$, where $f_{(ij)}$ are differentiable. It then holds that $\theta_{ij}=g(\theta_i)+h(\theta_j)$, for all $i,j\in\nodeset$ and some functions $g$ and $h$.
\end{thm}
Hence, the model of the above theorem leads to the $p_1$ model in \eqref{eq:betarepd}:
\begin{coro}\label{coro:p1}
A dyadic independent ERGM with $n$ nodal parameters obtained from a permutation-equivariant parametrization is a $p_1$ model.
\end{coro}

Similarly, for the directed dyadic-independent stochastic block models, we do not have the constraint $h_{(kl)}=h_{(lk)}$. further parametrization to $r<n$ block parameters is of form $\eta_{kl}=h_{(kl)}(\eta_k,\eta_l)$, $k,l\in[r]$.   Similar to the previous case, for a permutation $\pi$, we define $\pi(\eta_k,\eta_l)=(\eta_{\pi(k)},\eta_{\pi(l)})$; and similarly $\pi(\eta_{kl})=\eta_{\pi(k)\pi(l)}$. Doubly permutation equivariance is defined in the same manner as in the case of undirected dyadic-independent stochastic block model. 

Now we have an analogous result to Theorem \ref{thm:sbm}:
\begin{thm}\label{thm:sbmd}
If a dyadic-independent stochastic block model with $r<n$ block parameters is obtained from a permutation-equivariant parametrization $\eta_{kl}=h_{(kl)}(\eta_k,\eta_l)$, $k,l\in[r]$, where  $h_{(kl)}$ are differentiable, then $\eta_{kl}=g(\eta_k)+h(\eta_l)$, for all $k,l\in [r]$ and some functions $g$ and $h$.
\end{thm}
Hence the model of the above theorem leads to the directed additive stochastic block model in \eqref{eq:addsbmd}:
\begin{coro}
A dyadic independent directed stochastic block model with $r<n$ block parameters obtained from a permutation-equivariant parametrization is a directed additive stochastic block model.
\end{coro}

\section{Summary and discussion}\label{sec:sum}
We have shown that for dyadic independent ERGMs, and under the  assumption of permutation equivariance for the parametrization of dyadic parameters into nodal parameters, the only possible resulting statistical models  are some well-known ERGMs: in the case of undirected networks we obtain the $\beta$ model and in the case of directed networks we obtain the $p_1$ model. For the case of further parametrizing between-block parameters, as opposed to dyadic parameters, to fewer than $n$ block parameters, and under the same assumption of permutation equivariance, for between-block and block parameters, the resulting models are additive stochastic block models for undirected networks, and directed additive stochastic block models for directed networks.

These results provide a strong theoretical justification for the use of these models in statistical network analysis. Instead of strong modeling assumptions that could potentially be imposed by these models, it is evident that these models must hold under natural assumptions of exponential family, dyadic independence, and permutation equivariance.

To the best of our knowledge, permutation equivariance has not been used as an assumption for parametrization elsewhere. The concept of equivariance has indeed been used as a property of estimators  in statistics; see, e.g., \citet{leh05} and also \citet{kee10} for location equivariance. More recently, permutation equivariance has been widely used as a property of neural networks (see, e.g, \citet{vic21}) and transformers (see, e.g, \citet{lee19}) in machine learning. 

For further parametrization to a lower-dimensional parameter space where the original parameters are of similar type, and the new parameters are of the same type, permutation equivariance is a very natural assumption: if a new parameter is related to some original parameters in a particular way, one expects other parameters be related to their own relevant original parameters in the same fashion. For this reason, further studies on other statistical models where permutation-equivariance could be used, and also the study of different types of equivariance for parametrization in statistical models are  interesting subjects to pursue.

In particular, for the nodal parametrization studied in this work, it is quite reasonable to assume that dyadic parameters behave similarly relative to their corresponding nodal parameters. Hence, arguably, this assumption is not only natural but necessary for obtaining well-behaved statistical network models with nodal parameters. The same, of course,  can be said about block parameters. 

In the context of dyadic parametrization of random networks, the idea of equivariance seems quite relevant in general.  Indeed other types of equivariance parametrizations could be used for dyadic parameters. For example, an idea is that  soft intervention on the nodes of the network must satisfy some type of equivariance on the dyadic parameters. This could facilitate causal interpretation on random networks.

Another idea for future work stems from the fact that the parametric model implied by permutation equivariance has the same form as that of the Rasch model \citep{ras80}. Imposing (permutation) equivariance assumption on Rasch models is interesting on its own, but the consequences of using equivariance Rasch models for (potentially bipartite) random networks can have important applications.
%

\bibliographystyle{elsarticle-num-names} 
\bibliography{bib}






\end{document}